\documentclass[
a4paper,
]
{amsart}

\usepackage{pdfsync}
\usepackage{tikz}
\usepackage{amsmath,amssymb,amsthm}

\usepackage{lmodern}

\usepackage{hyperref}
\colorlet{linkcolour}{red!50!black}
\hypersetup{
pdftitle={Simplified a priori Estimate for the Time Periodic Burgers' Equation}, 
pdfauthor={Magnus Fontes and Olivier Verdier},
colorlinks=true,
linkcolor=linkcolour, citecolor=linkcolour, pagecolor= linkcolour,
pagebackref=true,
}

\newcommand{\R}{\mathbb{R}}
\newcommand{\Z}{\mathbb{Z}}

\newcommand{\dd}{\,\mathrm d}

\newcommand{\bracket}[2]{\left\langle #1,#2 \right\rangle}
\newcommand{\duality}[2]{\bracket{#1}{#2}}

\newcommand{\psL}[2]{\left(#1,#2\right)}
\newcommand{\normL}[1]{\abs{#1}}

\newcommand{\norm}[1]{\left\lVert#1\right\rVert}
\providecommand{\abs}[1]{\left\lvert#1\right\rvert}

\newcommand{\dt}[1]{#1_{\sqrt{t}}}

\newcommand{\dts}[1]{#1_{\sqrt{t}*}}
\newcommand{\Dx}[1]{#1_x}

\newcommand{\Dsop}[1]{\operatorname{D}^{#1}}
\newcommand{\Dssop}[1]{\operatorname{D}^{#1}_{*}}
\newcommand{\Dhalfop}{\Dsop{\frac 12}}

\newcommand{\Hilbop}{\operatorname{\mathcal{H}}}
\newcommand{\Burgop}{\operatorname{T}}
\newcommand{\Linop}{\operatorname{\mathcal{L}}}

\newcommand{\Hilbert}{H}
\newcommand{\DTH}{\mathcal{D}(\T,H)}
\newcommand{\dTH}{\mathcal{D}'(\T,H^*)}

\newcommand{\Sob}[1]{\mathrm{H}^{\left(#1\right)}}

\newcommand{\Sobsob}[2]{\mathrm{H}^{\left(#1\right)\left(#2\right)}}
\renewcommand{\H}{\Sob{\frac 12,1}_0}

\newcommand{\Hs}{\Sob{-\frac 12,-1}}
\newcommand{\Hss}{\Sob{0,-1}}
\newcommand{\Leb}[1]{\mathrm{L}^{#1}}
\newcommand{\Four}[1]{\hat{#1}}

\newcommand{\Hilb}[1]{\widetilde{#1}}

\newcommand{\Q}{Q}
\newcommand{\I}{I}
\newcommand{\T}{\mathbb{T}}
\newcommand{\TS}{\T\times \R}

\newcommand{\TxI}{\T\times\I}

\newcommand{\puls}{2\pi}
\newcommand{\e}{\mathrm{e}}
\newcommand{\im}{\mathrm{i}}
\newcommand{\expn}[1]{\e^{ \im \puls #1 t}}

\newcommand{\sgn}{\operatorname{sgn}}
\newcommand{\Id}{\operatorname{\mathrm{Id}}}

\newtheorem{thm}{Theorem}

\newtheorem{lma}{Lemma}[section]

\newtheorem{prop}{Proposition}[section]
\theoremstyle{remark}
\newtheorem{rem}{Remark}[section]
\theoremstyle{definition}

\newenvironment{pr}[1][\proofname]{
\begin{proof}[#1]\mbox{}}{\end{proof}
\normalsize}

\newcommand{\figcaption}[1]{\caption{\small#1}}

\title{Simplified a priori Estimate for the Time Periodic Burgers' Equation}

\author{Magnus Fontes \and Olivier Verdier}

\begin{document}

\begin{abstract}
We present here a  version of the
existence and uniqueness result of time periodic solutions to the viscous Burgers equation with  irregular forcing terms (with Sobolev regularity  $-1$ in space).
The key result here is an a priori estimate which is simpler than the previously treated case of forcing terms with regularity $-\frac{1}{2}$ in time.
\end{abstract}

\maketitle

\section*{Introduction}
The study of the Burgers' equation has a long history starting with the seminal papers by Burgers \cite{Burgers}, Cole \cite{Cole}  and Hopf \cite{Hopf} where the Cole-Hopf transformation was introduced. The Cole-Hopf transformation transforms the homogeneous Burgers' equation into the heat equation.

More recently there have been several articles dealing with the forced Burgers' equation:
\begin{equation*}\label{simple-Burger}
u_t - \nu u_{xx} + u u_x = f
\end{equation*}
The vast majority treats the initial value problem in time with homogeneous Dirichlet or periodic space boundary conditions (see for instance \cite{Kreiss}).

Only recently has the question of the time-periodic forced Burgers' equation been tackled (\cite{Moser,WeinanE,Sinai,Fokas}). In most cases \cite{Moser,WeinanE} the authors are chiefly interested in the inviscid limit (the limit when the viscosity $\nu$ tends to zero). 

The closest related work to ours is that of Jauslin, Kreiss and Moser \cite{Moser} in which the authors show existence and uniqueness of a space and time periodic solution of the Burgers' equation for a space and time periodic forcing term which is smooth.

\section{Definitions}

In this section we recall some well known facts and fix some general notations. 

We will be concerned will \emph{time-periodic} solutions, of the Burgers' equation so we will use the following notation for the one-dimensional torus $\T$:

\begin{equation*}
	\T = \R / \Z\\
\end{equation*}

\subsection{Fractional Derivatives}

Given a Hilbert space $H$, we will denote the space of test functions with values in $H$ by by $\DTH$. Its dual space will be denoted by $\dTH$:
\[\dTH = (\DTH)^*\]

For any positive real number $s$ we may define the fractional derivative of order $s$ in the following way on the space of  time-periodic distributions $\dTH$ :
\[\Dsop{s} u = \sum_{k\in\Z} (\puls \im k)^s u_k \expn{k} = \sum_{k\in\Z} \abs{\puls \im k}^{s} \e^{\im \sgn(k) s \frac{\pi}{2}}u_k \expn{k}\]
where we have used the principal branch of the logarithm. The \emph{sign function} is defined as follows:
\[ 
\sgn(k) := \begin{cases} \frac{k}{\abs{k}} & \text{if $k\neq 0$}\\
0 &\text{if $k=0$}\end{cases} 
\] 
For $s=0$ we define $\Dsop{0}=\Id$.
$\Dsop{1}$ coincides with the usual differentiation operator on $\dTH$. The familiar composition property also holds: $\Dsop{s} \circ \Dsop{t} = \Dsop{s+t}$ for any $t,s\geq 0$.

The \emph{adjoint operator} of $\Dsop{s}$ is defined by using the conjugate of the multiplier of $\Dsop{s}$:
$$\Dssop{s} u = \sum_{k\in\Z} \abs{\puls \im k}^{s} \e^{-\im \sgn(k) s \frac{\pi}{2}}u_k \expn{k}$$

$\Dsop{s}$ and $\Dssop{s}$ are adjoints in the sense that for any $u\in\dTH$ and $\varphi\in\DTH$:
$$\duality{\Dsop{s} u}{\varphi} = \duality{u}{\Dssop{s} \varphi}$$
and similarly:
$$\duality{\Dssop{s} u}{\varphi} = \duality{u}{\Dsop{s} \varphi}$$

\subsection{Hilbert Transform}

The \emph{Hilbert transform} $\Hilbop$ is defined using the multiplier $-\im\sgn(k)$. For $u\in\dTH$ let 
$$\Hilbop u = \sum_{k\in\Z}-\im\sgn(k)\,u_k\, \expn{k}$$
 Simple computations then give:
\[ \Dhalfop_* = \Dhalfop\circ\Hilbop = \Hilbop\circ\Dhalfop \]

Notice that if $\Hilbert$ is a function space then $\Hilbop$ maps real functions to real functions.

The following properties will be useful in the sequel.

\begin{equation}\label{dtshilbeq}
\forall u\in\Sob{\frac 12}(\T,\Hilbert) \quad \psL{\Dhalfop u}{\Dhalfop_* \Hilbop u}_{\Leb{2}(\T,\Hilbert)} = - \norm{\Dhalfop u}^2_{\Leb{2}(\T,\Hilbert)}
\end{equation}
\[
\forall u \in\Leb{2}(\T,H)\quad \Re\Bigl(\psL{u}{\Hilbop(u)}_{\Leb{2}(\T,H)}\Bigr) = 0
\]
where $\Re$ denotes the real part of the expression.

\subsection{Fractional Sobolev Spaces}

We define fractional Sobolev spaces in the following manner, for any $s\in\R$:
$$\Sob{s}(\T,\Hilbert) = \Bigl\{u \in \dTH;\quad \sum_{k\in\Z}\abs{1+k^2}^{s}\norm{u_k}_{\Hilbert}^2 <\infty\Bigr\}$$

Of course $\Sob{0}(\T,\Hilbert) = \Leb{2}(\T,\Hilbert)$. When $s\geq 0$ then for an $u\in\Leb{2}(\T,\Hilbert)$: $u\in\Sob{s}(\T,\Hilbert) \iff \Dsop{s} u \in \Leb{2}(\T,\Hilbert)$. Moreover $\Sob{s}(\T,\Hilbert)$ is then a Hilbert space with the following scalar product:
$$\psL{u}{v} := \psL{u}{v}_{\Leb{2}(\T,\Hilbert)} + \psL{\Dsop{s} u}{\Dsop{s} v}_{\Leb{2}(\T,\Hilbert)}$$

The following classical result holds: $\bigl(\Sob{s}(\T,\Hilbert)\bigr)^* = \Sob{-s}(\T,\Hilbert^*)$.

\subsection{Anisotropic Fractional Sobolev Spaces}

We will now describe some useful Sobolev spaces defined on the interval 
\[\I=(0,1)\]

We will also fix a non-negative real number $s$:
\[s\geq 0\]

Let $\Sob{s}(\I)$ denote the usual fractional Sobolev space of real-valued s-times differentiable functions on $\I$.  
$\Sob{s}_0(\I)$ is the closure of $\mathcal{D}(\I)$ in $\Sob{s}(\I)$. In that case we have $\bigl(\Sob{s}_0(\I)\bigr)^* = \Sob{-s}(\I)$.
We will also use the following notations, for $\alpha$, $\beta$ nonnegative real numbers:
$$\Sobsob{\alpha}{\beta}(\TxI) := \Sob{\alpha}(\T, \Sob{\beta}(\I))$$
and
$$\Sob{\alpha,\beta}(\TxI) := \Sobsob{\alpha}{0}(\TxI)\cap \Sobsob{0}{\beta}(\TxI)$$
We also introduce $\Sob{\alpha,\beta}_0(\TxI)$ as the closure of $\mathcal{D}(\TxI)$ in $\Sob{\alpha,\beta}(\TxI$). It is clear that $\Sob{\alpha,\beta}_0(\TxI) = \Sobsob{\alpha}{0}(\TxI)\cap\Leb{2}(\T,\Sob{\beta}_0(\I))$. 

In the following sections we will use the following notations. The main working space will be $\H(\TxI)$, for which we use the more concise notation:
\[\H := \H(\TxI) \]
More generally, we will drop the dependence of the spaces on the domain $\TxI$, so that $\Leb{p}$ stands for $\Leb{p}(\TxI)$, $\Sob{\alpha,\beta}$ stands for $\Sob{\alpha,\beta}(\TxI)$, and so on.

We will also use the following notations, for $u\in\H$:

\begin{equation}\label{unotations}
\forall u \in \H\qquad
\begin{aligned}
	\Hilb{u} &:= \Hilbop u\\
	\dts{u} &:= \Dhalfop_* u\\
\end{aligned}
\end{equation}

\section{Interpolation and regularity}
\label{secinterp}
If $s_k(\xi)$ is the Fourier transform $s_{k}(\xi) = \Four{u}(k,\xi)$ of a distribution $u$ defined on $\TS$, we have the following H\"o{}lder inequality for any $\theta\in [0,1]$:
\begin{multline*}
\int_{\R}\sum_{k \in  \Z} \abs{k}^{2\alpha(1-\theta)} \abs{\xi}^{2\beta\theta} \abs{s_{k}(\xi)}^2 \dd \xi \leq \\
 \left(\int_\R \sum_{k \in  \Z} \abs{k}^{2\alpha} \abs{s_{k}(\xi)}^2 \dd \xi \right)^{1-\theta}\left(\int_\R \sum_{k \in \Z} \abs{\xi}^{2\beta} \abs{s_{k}(\xi)}^2\dd \xi\right)^{\theta}
 \end{multline*}
From this H\"o{}lder inequality we deduce
\begin{equation*}
\Sob{\alpha,\beta}(\TS) \hookrightarrow \Sob{(1-\theta)\alpha}(\T, \Sob{\theta\beta}(\R))
\end{equation*}
So using an extension operator from $\Sob{\theta\beta}(\I)$ to $\Sob{\theta\beta}(\R)$ one can prove the corresponding inclusion:
\begin{equation*}\label{interpolationInclusions}
\Sob{\alpha,\beta}(\TxI) \hookrightarrow \Sobsob{(1-\theta)\alpha}{\theta\beta}(\TxI)\end{equation*}

For $\alpha = 1/2$ and $\beta = 1$ and $\theta = \frac 13$ we get:
$$\H(\TxI) \subset \Sob{\frac 1 2, 1}(\TxI) \subset \Sobsob{\frac{1}{3}}{\frac{1}{3}}(\TxI)$$
Then the vectorial Sobolev inequalities yield:
\begin{equation*}\label{HincludeLfour}  \H(\TxI) \subset \Sobsob{\frac{1}{3}}{\frac{1}{3}} (\TxI)\hookrightarrow \Leb{4}(\T, \Sob{\frac 1 3}(\I)) \hookrightarrow \Leb{4}(\T, \Leb{4}(\I)) = \Leb{4}(\TxI) \end{equation*}
Here the injection $\Sobsob{\frac{1}{3}}{\frac{1}{3}}(\TxI)\hookrightarrow\Leb{4}(\T,\Sob{\frac{1}{3}})$ is compact and thus the injection $\H(\TxI)\hookrightarrow\Leb{4}(\TxI)$ is \emph{compact}.

  \tikzstyle{axes}=[]
  \tikzstyle{accessory}=[dashed]
  \tikzstyle{principal}=[thick]
  \tikzstyle{secondary}=[thick, dashed]
	\tikzstyle{intersection}=[black]

\begin{figure}
\begin{center}
\begin{tikzpicture}[scale=3]
  \begin{scope}[style=axes]
    \draw[->] (0,0) -- (.8,0) node[right] {$\alpha$};
    \draw[->] (0,0) -- (0,1.3) node[above] {$\beta$};
	\end{scope}
    \foreach \x/\xtext in {.3333333333333/\frac{1}{3}, .5/\frac{1}{2}}
      \draw[xshift=\x cm] (0pt,1pt) -- (0pt,-1pt) node[below,fill=white]
            {$\xtext$};
    \foreach \y/\ytext in {.3333333333333/\frac{1}{3}, 1}
      \draw[yshift=\y cm] (1pt,0pt) -- (-1pt,0pt) node[left,fill=white]
            {$\ytext$};
	\draw[principal] (0,1) -- (.5,0);
	\draw[accessory] (0,0) coordinate (orig) -- (.5,.5) coordinate (equil);
	\fill[red] (intersection of orig--equil and .5,0--0,1) circle(.5pt);
\end{tikzpicture}
\figcaption{$\H$ is included in $\Sobsob{\frac{1}{3}}{\frac{1}{3}}$ which is included in $\Leb{6}$ by the usual Sobolev inclusion theorem. In particular, $\H$ is included in $\Leb{4}$, so $u\in\H\implies u^2\in\Leb{2}$. As a result the non-linear term of the Burgers' equation may be written as $-(u^2,v_x)$ for a test function $v\in\H$ since $v\in\H\implies v_x\in\Leb{2}$ by definition.}\label{figLfourregularity}
\end{center}
\end{figure}
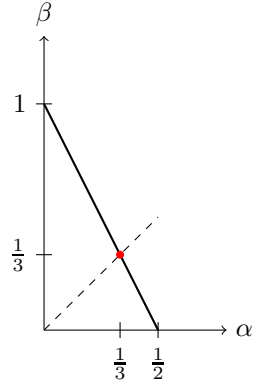

\section{Main Result: a simplified a priori estimate}

We define the Burgers' Operator by:
\[
\Burgop = \Linop + S
\]
where $\Linop$ and $S$ are defined in the familiar weak form, the bracket being the \emph{duality bracket} between $\H$ and $\Hs$ (recall the notations \eqref{unotations}):
\[\forall v \in \H \quad \duality{\Linop u}{v} := \psL{\dt u}{\dts v} + \mu \psL{u_x}{v_x}
\]
and
\[\forall v \in \H\quad \duality{S(u)}{v} := -\frac 12 \psL{{u^2}}{\Dx{v}}\]
It turns out that the second definition makes sense because of the embedding $\H\subset \Leb{4}$ (see \autoref{figLfourregularity}).



\newcommand{\NTS}{N_{\TS}}
\newcommand{\homog}{\mathcal{C}}
\newcommand{\NQ}{N_Q}
\newcommand{\Rez}{R_0}

We now prove the main result of this article, namely an a priori estimate on the solutions of the family of equations:
\[(\Linop + \lambda S)u = f\]
for $f\in \Hss$ and $0\leq \lambda \leq 1$.

In the proof we use techniques similar to those in \cite{Fontes-Saksman}. The proof is simpler (but the result is weaker) than the one obtained in \cite{fontes-verdier}.

\begin{thm}\label{lambdaaprioriestimatethm}
Let $f\in\Hss$. The set 
$$\bigcup_{\lambda\in [0,1]} (\Linop +\lambda S)^{-1}\bigl(\{f\}\bigr)$$
is bounded in $\H$.
\end{thm}

We will need the following Lemma which may be proved using a scaling argument (a proof is available in \cite{fontes-verdier}).
\begin{lma}\label{ThmHomog}
There exists a constant $\homog\in\R$ such that for any $u\in \H(\Q)$:
$$\int_\Q \abs{u(t,x)}^4 \dd t\dd x \leq \homog^2 \left(\int_\Q \abs{u}^2 \dd t \dd x + \int_\Q \abs{\dt{u}}^2\dd t\dd x\right) \cdot \left(\int_\Q \abs{u_x}^2\dd t\dd x \right)$$
which implies that:
\begin{equation*}\label{usqlequux}
\normL{u^2}\leq \homog \norm{u} \normL{u_x}
\end{equation*}
\end{lma}

\begin{pr}[Proof of \autoref{lambdaaprioriestimatethm}]
By definition $\Linop u + \lambda S(u) = f$ means:
\begin{equation}
\forall v\in \H\quad \psL{\dt{u}}{\dts{v}} + \mu\psL{u_x}{v_x} -\frac 12 \lambda\psL{u^2}{v_x} = \langle f, v \rangle
\label{burgersdefvareq}
\end{equation}
\begin{enumerate}
\item
We notice that for smooth $u$:
\begin{align*}
\psL{u^2}{u_x} &= \int_{\Q} u^2 u_x\\
				&= \frac{1}{3}\int_{\Q} (u^3)_x \\
				&= 0
\end{align*}
and then by density and continuity this holds for all $u\in\H$.
\item
With $v=u$ in \eqref{burgersdefvareq} we get:
$$\underbrace{\psL{\dt{u}}{\dts{u}}}_{=0} + \mu \psL{u_x}{u_x} + \frac{1}{2}\lambda\underbrace{\psL{u^2}{u_x}}_{=0} = \langle f, u \rangle$$
which gives:
\[\normL{u_x}^2 = \frac{\langle f, u \rangle}{\mu }\label{muuxequalf}\\
			\leq \frac{\norm{f}\normL{u_x}}{\mu}\]
From this we deduce that
\begin{equation}\label{gooduxestimate}
\normL{u_x} \leq \frac{\norm{f}}{\mu}
\end{equation}
\item
Pairing in \eqref{burgersdefvareq} with the Hilbert transform of $u$, $v=\Hilb{u}$  we get:
\[\psL{\dt{u}}{\dts{\Hilb u}} + \mu \underbrace{\psL{u_x}{\Hilb{u}_x}}_{=0} + \frac{1}{2} \lambda \psL{u^2}{\Hilb{u}_x} = \langle f, \Hilb{u} \rangle\]
Using the identity \eqref{dtshilbeq}, the fact that $\normL{\Hilb{u_x}} = \normL{u_x}$ and that $\lambda \leq 1$ we get:
\begin{equation*}
\normL{\dt{u}}^2\leq \frac{1}{2}\abs{\psL{u^2}{\Hilb{u}_x}} + \norm{f}\normL{u_x}
\label{dtuestimate}
\end{equation*}
\item
We estimate $\abs{\psL{u^2}{{\Hilb{u}}_x}} $ using \autoref{ThmHomog}:
\begin{equation}\begin{split}
\abs{\psL{u^2}{{\Hilb{u}}_x}} &\leq  \normL{u^2}\normL{u_x} \\
				&\leq  \homog \norm{{u}} \normL{u_x}^2 \label{rawconvectionestimate}
\end{split}\end{equation}
\item
Using the estimate \eqref{gooduxestimate} inside \eqref{rawconvectionestimate} we obtain:
\begin{equation*}\label{dtuRoestimate}
\begin{split}
\normL{\dt{u}}^2 &\leq  \frac{\homog}{2}\norm{f}\normL{u_x}^2 + \norm{f}\normL{u_x}\\
&\leq \frac{\norm{f}^2}{\mu}\biggl(\frac{\homog}{2\mu}\norm{u}+1\biggr)
\end{split}
\end{equation*}

Since that estimate does not depend on $\lambda$ the theorem is proved.
\end{enumerate}
\end{pr}

The a priori estimate above may now be used to prove existence of solutions by a (nonlinear, compact) degree argument using the Leray-Schauder Theorem (cf. \cite{fontes-verdier}).

\begin{figure}
\begin{center}
\begin{tikzpicture}[scale=3]
  \begin{scope}[style=axes]
    \draw[->] (0,0) -- (1.3,0) node[right] {$\alpha$};
    \draw[->] (0,0) -- (0,2.3) node[above] {$\beta$};
	\end{scope}
    \foreach \x/\xtext in {.6666666666666/\frac{2}{3}, .5/\frac{1}{2}, 1}
      \draw[xshift=\x cm] (0pt,1pt) -- (0pt,-1pt) node[below,fill=white]
            {$\xtext$};
    \foreach \y/\ytext in {.6666666666666/\frac{2}{3},  1, 2}
      \draw[yshift=\y cm] (1pt,0pt) -- (-1pt,0pt) node[left,fill=white]
            {$\ytext$};
	\draw[principal] (0,2) -- (.5,1);
	\draw[secondary] (.5,1) -- (1,0);
	\draw[accessory] (0,0) coordinate (orig) -- (.8,.8) coordinate (equil);
	\fill[red] (intersection of orig--equil and 1,0--0,2) circle(.5pt);
\end{tikzpicture}
\figcaption{The first step of the Cole-Hopf Transformation is an integration in $x$. This function $U$ obtained thus ends up in $\Sobsob{0}{1}\cap\Sobsob{\frac{1}{2}}{1}$, which delimits the plain line on the graph above. But it follows from $\Burgop u \in \Sobsob{0}{-1}$ that $u$ is actually also in $\Sobsob{1}{-1}$ so $U$ ends up in $\Sobsob{1}{2}$ and we have an inclusion in $\Sobsob{\frac{2}{3}}{\frac{2}{3}}$ which is embedded in continuous H\"older functions.} 
\label{colehopf}
\end{center}
\end{figure}
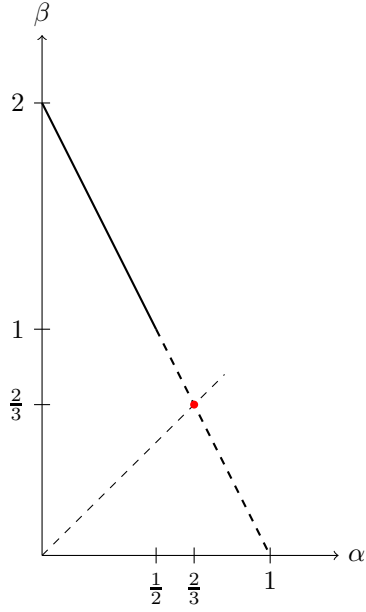

\section{Existence and Uniqueness of the time-periodic solutions}

A weaker result of the main result proved in \cite{fontes-verdier} is
\begin{thm}
For $f\in\Sobsob{0}{-1}$ there exists a unique solution $u\in\H$ of
\[\Burgop u = f\]
\end{thm}

In this section we briefly sketch the proof of that result, the full proof being available in \cite{fontes-verdier}.

The main ingredient of the proof is the Cole-Hopf transformation, which is essentially defined by:
\[u = \frac{\varphi_x}{\varphi}\]

In our case there are complications due to the fact that $u\in\H$, so $u$ is time-periodic. This change of variable will transform the periodicity problem into an eigenvalue problem (because the Cole-Hopf transformation linearises the Burgers' equation). After working out the details one shows that the uniqueness problem is equivalent to the uniqueness of the \emph{ground state eigenvalue problem}:

\begin{prop}
Given $v\in\H$ the solution set of the following equation in $K$ and $\varphi$
\begin{equation*}\begin{cases}
\varphi_t -\mu\varphi_{xx} + v\varphi_x + K\varphi = 0\\
\varphi >0\\
\varphi_x|_{\partial \Q} = 0\\
\varphi \in\Sob{1,2}\\
K\in\R
\end{cases}\end{equation*}
is $K=0$ and $\varphi = 1$ \emph{if and only if} $\Burgop u = \Burgop v$ implies $u = v$ (that is, the solution to the original Burger's equation is unique).
\end{prop}

The proof of that proposition essentially hinges on the embedding properties exposed in \autoref{secinterp} (see \autoref{colehopf}).

The remaining part of the proof is concerned with the eigenvalue problem of the Proposition above. One first shows that the eigenvalue is zero using a weaker version of the Perron-Frobenius theorem. The second step is to show that the remaining eigenvalue problem is \emph{non degenerate}, namely that the dimension of the eigenspace must be one. This last step makes use of the a priori estimate proved in \autoref{lambdaaprioriestimatethm}.

The details of that part of the proof are too lengthy to be exposed here in depth so the interested reader is referred to \cite{fontes-verdier}.


\end{document}